\documentclass[12pt]{article}

\usepackage{amstext}
\usepackage{latexsym}
\usepackage{amssymb}

\def\R{\mathbb R}

\def\N{\mathbb N}

\newcommand{\eps}{\varepsilon}

\begin{document}

\newtheorem{theorem}{Theorem}[section]
\renewcommand{\thetheorem}{\arabic{section}.\arabic{theorem}}
\newtheorem{definition}[theorem]{Definition}
\newtheorem{deflem}[theorem]{Definition and Lemma}
\newtheorem{lemma}[theorem]{Lemma}
\newtheorem{example}[theorem]{Example}
\newtheorem{remark}[theorem]{Remark}
\newtheorem{remarks}[theorem]{Remarks}
\newtheorem{cor}[theorem]{Corollary}
\newtheorem{pro}[theorem]{Proposition}
\newtheorem{proposition}[theorem]{Proposition}

\renewcommand{\theequation}{\thesection.\arabic{equation}}

\title{Higher regularity of the ``tangential'' fields 
\\ in the relativistic Vlasov-Maxwell system}

\date{}

\author{{\sc Markus Kunze}\\[2ex]
         Mathematisches Institut, Universit\"at K\"oln, \\
         Weyertal 86-90, D\,-\,50931 K\"oln, Germany \\
         e-mail: mkunze@mi.uni-koeln.de}

\maketitle

\begin{abstract}\noindent
It is shown that the ``tangential'' electric and magnetic fields, 
in the Glassey-Strauss representation formulas, are in fact 
bounded in $L^\infty_{{\rm loc},\,t} L^{2+\delta}_x$ for some $\delta>0$. 
\end{abstract}


\setcounter{equation}{0}

\section{Introduction and main result}

The relativistic Vlasov-Maxwell system describes the time evolution 
of a plasma with particles moving at high velocities 
(close to the speed of light which is taken to be $c=1$).  
The Vlasov equation
\begin{equation}\label{vlas1}
   \partial_t f+v\cdot\nabla f+(E+v\wedge B)\cdot\nabla_p f=0
\end{equation}
governs the evolution of the scalar density function $f=f(t, x, p)\ge 0$, 
depending on time $t\in\R$, position $x\in\R^3$, and momentum $p\in\R^3$; 
here $\nabla$ always means $\nabla_x$.
The velocity $v\in\R^3$ associated to $p$ is
\[ v=\frac{p}{\sqrt{1+p^2}}\,,\quad\mbox{thus}\quad p=\frac{v}{\sqrt{1-v^2}}\,, \] 
where $p^2=|p|^2$ and $v^2=|v|^2$ for brevity. The Lorentz force
\[ L=L(t, x, v)=E(t, x)+v\wedge B(t, x)\in\R^3 \]
is obtained from the electric field $E=E(t, x)\in\R^3$
and the magnetic field $B=B(t, x)\in\R^3$, which in turn satisfy
the Maxwell equations
\begin{equation}\label{maxE}
   \partial_t E=\nabla\wedge B-j,\quad\nabla\cdot E=\rho,
\end{equation}
and
\begin{equation}\label{maxB}
   \partial_t B=-\nabla\wedge E,\quad\nabla\cdot B=0.
\end{equation}
The coupling of (\ref{vlas1}) to (\ref{maxE}), (\ref{maxB}) is realized
through the charge density $\rho=\rho(t, x)\in\R$ and the current density
$j=j(t, x)\in\R^3$ via
\[ \rho(t, x)=\int_{\R^3}f(t, x, p)\,dp\quad\mbox{and}\quad
   j(t, x)=\int_{\R^3} v\,f(t, x, p)\,dp. \] 
Furthermore, initial data
\[ f(t=0)=f^{(0)},\quad E(t=0)=E^{(0)},\quad\mbox{and}\quad
   B(t=0)=B^ {(0)} \] 
are prescribed such that the constraint equations
\[ \nabla\cdot E^{(0)}=\rho^{(0)}=\int_{\R^3} f^{(0)}\,dp
   \quad\mbox{and}\quad\nabla\cdot B^{(0)}=0 \] 
are satisfied.
\smallskip 

There has been quite some activity concerning the relativistic Vlasov-Maxwell 
over the years, but nonetheless the question whether (for instance smooth) 
initial data will yield a global in time solution still remains open. 
See \cite{glassey} and \cite{strauss} for a general introduction and overview, 
\cite{Kyet} for a summary of results up to approximately 2015 
and \cite{patel} for some newer and further refined criteria 
concerning unrestricted global existence. 
\smallskip 

To explain the observation which is the subject of the present paper 
first recall that the energy
\begin{equation}\label{energcons} 
   {\cal E}(t)=\int_{\R^3}\int_{\R^3}\sqrt{1+p^2}\,f(t, x, p)\,dx\,dp
   +\frac{1}{2}\int_{\R^3}(|E(t, x)|^2+|B(t, x)|^2)\,dx
\end{equation}  
is conserved along solutions of (\ref{vlas1}), (\ref{maxE}), and (\ref{maxB});
note that $\nabla_p\sqrt{1+p^2}=v$. Therefore one gets a bound
\begin{equation}\label{energybd}
   E, B\in L^\infty_t L^2_x
\end{equation} 
in terms of the initial data for free. Next, defining 
$E^{(1)}(x)=\partial_t E(0, x)$ and $B^{(1)}(x)=\partial_t B(0, x)$, 
$E$ and $B$ are the solutions to the wave equations
\begin{eqnarray}
   \square\,E & = & -(\partial_t j+\nabla\rho)
   =-\int_{\R^3} (v\,\partial_t+\nabla)f\,dp,
   \quad E(0)=E^{(0)},\quad \partial_t E(0)=E^{(1)},
   \label{E-wave}
   \\ \square\,B & = & \nabla\wedge j
   =\nabla\wedge\int_{\R^3} v\,f\,dp,
   \hspace{6.5em} B(0)=B^{(0)},\quad \partial_t B(0)=B^{(1)}.
   \label{B-wave}
\end{eqnarray}
In the pioneering paper \cite{glstr}, Glassey and Strauss 
noted that (\ref{E-wave}) and (\ref{B-wave}) can be used to derive 
representation formulas for the fields as follows. Write 
\[ S=\partial_t+v\cdot\nabla,\quad T_j=-\omega_j\partial_t+\partial_{x_j}. \]   
Then $\partial_t$ and $\nabla$ can be expressed in terms of $S$ and $T$, since 
\begin{eqnarray} 
   \partial_t & = & (1+v\cdot\omega)^{-1} (S-v\cdot T),
   \label{partt} 
   \\ \partial_{x_j} & = & T_j+(1+v\cdot\omega)^{-1}\omega_j (S-v\cdot T). 
   \label{partxj} 
\end{eqnarray}
Note that with $\omega=\frac{y-x}{|y-x|}$: 
\[ \nabla_y\,[f(t-|y-x|, y, p)]=(-\omega\,\partial_t+\nabla)f(\ldots)
   =(Tf)(\ldots),  \] 
and, for instance,  
\begin{eqnarray*} 
   E & = & -\,\square^{-1}\int_{\R^3} dv\,(\nabla+v\,\partial_t)f
   \\ & \cong & -\int_{\R^3} dv\,\int_{|y-x|\le t}\frac{dy}{|y-x|}
   \,(\nabla+v\,\partial_t)f.
\end{eqnarray*} 
First one uses (\ref{partxj}) and (\ref{partt}) for the right-hand side 
and then one integrates $(Tf)(\ldots)=\nabla_y\,[\ldots]$ by parts in $y$. 
After a lengthy calculation one finds 
\begin{equation}\label{GlStrrep} 
   E = E_D+E_{DT}+E_T+E_S
\end{equation}  
(and a similar expression for $B$), where $E_D$ and $E_{DT}$ are data terms, 
\begin{eqnarray}
   E_T(t, x) & = & -\int_{|y|\le t}\frac{dy}{|y|^2}
   \int_{\R^3} dp\,K_{E,\,T}(\omega, v) f(t-|y|, x+y, p),
   \label{ETform}  
   \\ E_S(t, x) & = & -\int_{|y|\le t}\frac{dy}{|y|}
   \int_{\R^3} dp\,K_{E,\,S}(\omega, v)\,(Lf)(t-|y|, x+y, p),
   \label{ESform} 
\end{eqnarray}
and the integral kernels $K_{E,\,T}(\omega, v)\in\R^3$ 
and $K_{E,\,S}(\omega, v)\in\R^{3\times 3}$ behave as follows:  
\begin{eqnarray}
   & & |K_{E,\,T}(\omega, v)|\le C(1+p^2)^{-1}(1+v\cdot\omega)^{-3/2},
   \label{beh1} \\[1ex] 
   & & |K_{E,\,S}(\omega, v)z|
   \le C(1+p^2)^{-1/2}(1+v\cdot\omega)^{-1}|z|
   \quad (z\in\R^3).
   \label{beh2}
\end{eqnarray}
See Section \ref{kersect} below for the precise form of the kernels 
and a recap of the proof of (\ref{beh1}) and (\ref{beh2}). 
Relation (\ref{GlStrrep}) is the Glassey-Strauss representation formula 
for the electric field $E$ and, together with its counterpart for $B$, 
it has become an indipensible tool for proving existence results 
for the relativistic Vlasov-Maxwell system. Variants of it 
have been used for related systems as well. 

Assuming initial data of compact support, certainly the data terms 
$E_D$ and $E_{DT}$ in (\ref{GlStrrep}) will behave well. 
Thus, in the light of (\ref{energybd}), it is natural to ask 
what could be said about the terms $E_T$ and $E_S$ individually. 
We will call $E_T$ the tangential part and we are going to prove 
the following result. 

\begin{theorem}\label{mainthm} 
Consider initial data of compact support. 
Then $E_T, B_T\in L^\infty_{{\rm loc}, t} L^{2+\delta}_x$ 
for some $\delta>0$. 
\end{theorem} 

\begin{remark}{\rm (a) Since the argument for $B_T$ is the same 
as for $E_T$, we will only consider the latter in what follows. 
\smallskip 

\noindent 
(b) The number $\delta>0$ will be a uniform 
constant, for instance $\delta=\frac{2}{17}$ is a possible choice. 
As this result is mainly understood to be a ``proof of concept'', 
certainly the regularity that is gained here will not be optimal. 
\smallskip 

\noindent 
(c) By $E_T\in L^\infty_{{\rm loc}, t} L^{2+\delta}_x$ we mean 
the following: There is a continuous function $C=C(t): [0, \infty[\to [1, \infty[$ 
which only depends on $t$ and the initial energy ${\cal E}(0)$, 
the initial mass ${\cal M}(0)=\int_{\R^3}\int_{\R^3}f^{(0)}(x, p)\,dx\,dp$ 
and ${\|f^{(0)}\|}_\infty$ such that ${\|E_T(t, \cdot)\|}_{L^{2+\delta}_x(\R^3)}
\le C(t)$ for $t\in [0, T_{{\rm max}}[$, where $T_{{\rm max}}>0$ 
denotes the maximal time of existence of the solution. 
A constant denoted by $C$ will always be one which only depends on 
${\cal E}(0)$, ${\cal M}(0)$ and ${\|f^{(0)}\|}_\infty$.  
\smallskip 

\noindent 
(d) Due to Theorem \ref{mainthm} and (\ref{energybd}) one has 
$E_S\in L^\infty_{{\rm loc}, t} L^2_x$, but we are not able 
to derive this bound directly from (\ref{ESform}). 
}
\end{remark}


\setcounter{equation}{0}

\section{Proof of Theorem \ref{mainthm}} 

According to (\ref{ETform}) and (\ref{beh1}) we have 
\begin{eqnarray*} 
   |E_T(t, x)| & \le & 
   C\int_{|y|\le t}\frac{dy}{|y|^2}\int_{\R^3}\frac{dp}{1+p^2}
   \,\frac{1}{(1+v\cdot\omega)^{3/2}}\,f(t-|y|, x+y, p)
   \\ & =: & Cu(t, x). 
\end{eqnarray*} 
The Fourier transform of $u$ is 
\begin{eqnarray*} 
   \hat{u}(t, \xi) & = & \int_{\R^3} e^{-i\,\xi\cdot x}\,u(t, x)\,dx
   \\ & = & \int_{|y|\le t}\frac{dy}{|y|^2}\int_{\R^3}\frac{dp}{1+p^2}\,\frac{1}{(1+v\cdot\omega)^{3/2}}
   \,\int_{\R^3} dx\,e^{-i\,\xi\cdot x}\,f(t-|y|, x+y, p)
   \\ & = & \int_{\R^3}\frac{dp}{1+p^2}\,\int_{|y|\le t}\frac{dy}{|y|^2}\,e^{i\,\xi\cdot y}
   \,\frac{1}{(1+v\cdot\omega)^{3/2}}\,\hat{f}(t-|y|, \xi, p)
   \\ & = & \int_{\R^3}\frac{dp}{1+p^2}\,\int_0^t ds\,\hat{f}(t-s, \xi, p)
   \,\int_{|\omega|=1} dS(\omega)\,\frac{e^{is\,\xi\cdot\,\omega}}{(1+v\cdot\omega)^{3/2}}. 
\end{eqnarray*} 
To evaluate the inner integral choose a unit vector $e\in\R^3$ such that $\{\bar{v}, \bar{u}, e\}$ 
is an orthonormal basis of $\R^3$, where $\bar{v}=v/|v|=p/|p|$ 
and $\bar{u}=\frac{\bar{\xi}-(\bar{\xi}\cdot\bar{v})\bar{v}}{\sqrt{1-(\bar{\xi}\cdot\bar{v})^2}}$ 
are orthogonal unit vectors. Consider the matrix 
$A=\left(\begin{array}{c} e \\ \bar{u} \\ \bar{v}\end{array}\right)\in\R^{3\times 3}$, 
where the vectors are taken as rows. Then $A\bar{v}=e_3$ and $A\bar{u}=e_2$. 
It follows that $Av=|v|A\bar{v}=|v|e_3$ 
and $A\xi=|\xi|A\bar{\xi}=|\xi|A(\sqrt{1-(\bar{\xi}\cdot\bar{v})^2}\,\bar{u}+(\bar{\xi}\cdot\bar{v})\bar{v})
=|\xi|(\sqrt{1-(\bar{\xi}\cdot\bar{v})^2}\,e_2+(\bar{\xi}\cdot\bar{v})e_3)$, which in turn yields 
\begin{eqnarray*} 
   \int_{|\omega|=1} dS(\omega)\,\frac{e^{is\,\xi\cdot\,\omega}}{(1+v\cdot\omega)^{3/2}}
   & = & \int_{|\omega|=1} dS(\omega)\,\frac{e^{is\,|\xi|(\sqrt{1-(\bar{\xi}\cdot\bar{v})^2}\, 
   \omega_2+(\bar{\xi}\cdot\bar{v})\omega_3)}}{(1+|v|\omega_3)^{3/2}}
   \\ & = & \int_0^{2\pi} d\theta\int_0^{\pi} d\varphi\,\sin\varphi
   \,\frac{e^{is\,|\xi|(\sqrt{1-(\bar{\xi}\cdot\bar{v})^2}\,\sin\theta\sin\varphi
   +(\bar{\xi}\cdot\bar{v})\cos\varphi)}}
   {(1+|v| \cos\varphi)^{3/2}}
   \\ & = & \int_{-1}^1 d\sigma\,\frac{e^{is\sigma\,\xi\cdot\bar{v}}}{(1+|v|\sigma)^{3/2}}
   \int_0^{2\pi} d\theta\,e^{is |\xi|\sqrt{1-(\bar{\xi}\cdot\bar{v})^2}\,\sin\theta\sqrt{1-\sigma^2}}
   \\ & = & 2\pi\int_{-1}^1 d\sigma\,\frac{e^{is\sigma\,\xi\cdot\bar{v}}}{(1+|v|\sigma)^{3/2}}
   J_0\Big(s |\xi|\sqrt{1-(\bar{\xi}\cdot\bar{v})^2}\,\sqrt{1-\sigma^2}\Big), 
\end{eqnarray*} 
where 
\[ J_0(r)=\frac{1}{2\pi}\int_0^{2\pi} e^{ir\sin\theta}\,d\theta \] 
is the Bessel function of order zero. Its asymptotic expansion is 
\begin{equation}\label{bessel-expa} 
   J_0(r)=\sqrt{\frac{2}{\pi r}}\cos\Big(r-\frac{\pi}{4}\Big)+O(r^{-3/2}),\quad r\to\infty,
\end{equation} 
see \cite[p.~432]{grafa}, and also $|J_0(r)|\le 1$ is verified. Thus altogether we obtain 
\[ \hat{u}(t, \xi)=2\pi\int_0^t ds\int_{\R^3}\frac{dp}{1+p^2}\,\hat{f}(t-s, \xi, p)
   \,\int_{-1}^1 d\sigma\,\frac{e^{is\sigma\,\xi\cdot\bar{v}}}{(1+|v|\sigma)^{3/2}}
   J_0\Big(s |\xi|\sqrt{1-(\bar{\xi}\cdot\bar{v})^2}\,\sqrt{1-\sigma^2}\Big). \]
Next we introduce a standard Littlewood-Paley decomposition of $u$. 
For, fix $\varphi_0\in C_0^\infty(\R^n)$ such that $\varphi_0(\xi)=1$
for $|\xi|\le 1$ and $\varphi_0(\xi)=0$ for $|\xi|\ge 2$.
For $j\in\N$ put $\varphi_j(\xi)=\varphi_0(2^{-j}\xi)-\varphi_0(2^{-j+1}\xi)$.
Then $\varphi_j(\xi)=0$ for $|\xi|\le 2^{j-1}$ and for $|\xi|\ge 2^{j+1}$.
Furthermore, $\sum_{j=0}^\infty\varphi_j(\xi)=1$ for all $\xi\in\R^n$. 
Henceforth we shall consider $u_j=u_j(t, x)$ given by 
$\hat{u}_j(t, \xi)=\varphi_j(\xi)\hat{u}(t, \xi)$ for $j\in\N_0$. 
In this way we obtain 
\[ u=\sum_{j=0}^\infty u_j \] 
for 
\[ \hat{u}_j(t, \xi)=2\pi\int_0^t ds\int_{\R^3}\frac{dp}{1+p^2}\,\hat{f}_j(t-s, \xi, p)
   \,\int_{-1}^1 d\sigma\,\frac{e^{is\sigma\,\xi\cdot\bar{v}}}{(1+|v|\sigma)^{3/2}}
   J_0\Big(s |\xi|\sqrt{1-(\bar{\xi}\cdot\bar{v})^2}\,\sqrt{1-\sigma^2}\Big), \]
where $\hat{f}_j(t, \xi, p)=\varphi_j(\xi)\hat{f}(t, \xi, p)$; 
the Fourier transform of $f$ only refers to the variable $x$. Then 
\begin{equation}\label{qmon}
   {\|f_j(t, \cdot, p)\|}_{L^q_x(\R^3)}\le C {\|f(t, \cdot, p)\|}_{L^q_x(\R^3)},
   \quad j\in\N_0,\quad q\in [1, \infty], 
\end{equation}
uniformly in $t$ and $p$; the constant $C>0$ does only depend on $q$. 
Since ${\rm supp}\hat{f}_j(t, \cdot, p)\subset\{2^{j-1}\le |\xi|\le 2^{j+1}\}$,
Bernstein's inequality (or a direct estimate) moreover leads to 
\begin{equation}\label{bern1}
   {\|f_j(t, \cdot, p)\|}_{L^2_x(\R^3)}\le C 2^{3j/2}{\|f_j(t, \cdot, p)\|}_{L^1_x(\R^3)}, 
\end{equation} 
uniformly in $t$ and $p$. Denote by ${(\psi_j)}_{j\in\N_0}$ 
a partition of unity on $]0, 1]$ such that ${\rm supp}\,\psi_0 \subset [\frac{1}{3}, 1]$ 
and ${\rm supp}\,\psi_j\subset [2^{-(j+2)}, 2^{-j+1}]$
for $j\in\N$. Accordingly we decompose 
\begin{equation}\label{uhat-dec} 
   \hat{u}_j(t, \xi)=\sum_{k=0}^\infty\sum_{m=0}^\infty
   \sum_{n=0}^\infty \hat{u}_{jkmn}(t, \xi),
\end{equation} 
where
\begin{eqnarray*} 
   \hat{u}_{jkmn}(t, \xi) & = & 2\pi\int_0^t ds\,\psi_k\Big(\frac{s}{t}\Big)\int_{\R^3}\frac{dp}{1+p^2}
   \,\hat{f}_j(t-s, \xi, p)\,\psi_n\Big(\sqrt{1-(\bar{\xi}\cdot\bar{v})^2}\Big)
   \\ & & \hspace{2.5em} \times\,\int_{-1}^1 d\sigma\,
   \frac{e^{is\sigma\,\xi\cdot\bar{v}}}{(1+|v|\sigma)^{3/2}}
   J_0\Big(s |\xi|\sqrt{1-(\bar{\xi}\cdot\bar{v})^2}\,\sqrt{1-\sigma^2}\Big)\,\psi_m(\sqrt{1-\sigma^2}). 
\end{eqnarray*} 

The next lemma is the main technical tool 
for the proof of Theorem \ref{mainthm}. 

\begin{lemma}\label{dreifrz} For $j\in\N$ and $k, m, n\in\N_0$, 
\begin{equation}\label{ujkmnL2} 
   {\|\hat{u}_{jkmn}(t, \cdot)\|}_{L^2_\xi(\R^3)}
   \le Ct\min\Big\{1, \frac{2^{(k+m+n-j)/2}}{t^{1/2}}\Big\}\,2^{-k}
   \min\Big\{2^{-2m}\,2^{3j/2}, (\sqrt{n}+\sqrt{j})\,2^{-n}\Big\}.  
\end{equation} 
\end{lemma} 
{\bf Proof\,:} Observe that by (\ref{bessel-expa}) always 
\begin{eqnarray*} 
   \lefteqn{\psi_k\Big(\frac{s}{t}\Big)
   \psi_n\Big(\sqrt{1-(\bar{\xi}\cdot\bar{v})^2}\Big)\psi_m(\sqrt{1-\sigma^2})
   \,\Big|J_0\Big(s |\xi|\sqrt{1-(\bar{\xi}\cdot\bar{v})^2}\,\sqrt{1-\sigma^2}\Big)\Big|}
   \\ & \le & C\psi_k\Big(\frac{s}{t}\Big)\psi_n\Big(\sqrt{1-(\bar{\xi}\cdot\bar{v})^2}\Big)
   \psi_m(\sqrt{1-\sigma^2})
   \,\min\Big\{1, \frac{1}{s^{1/2} |\xi|^{1/2} 
   (1-(\bar{\xi}\cdot\bar{v})^2)^{1/4}\,(1-\sigma^2)^{1/4}}\Big\}
   \\ & \le & C\psi_k\Big(\frac{s}{t}\Big)\psi_n\Big(\sqrt{1-(\bar{\xi}\cdot\bar{v})^2}\Big)
   \psi_m(\sqrt{1-\sigma^2})
   \,\min\Big\{1, \frac{2^{(k+m+n-j)/2}}{t^{1/2}}\Big\}. 
\end{eqnarray*} 
Therefore 
\begin{eqnarray}\label{hatu-kmn} 
   \lefteqn{|\hat{u}_{jkmn}(t, \xi)|} 
   \nonumber \\ & \le & C\min\Big\{1, \frac{2^{(k+m+n-j)/2}}{t^{1/2}}\Big\}
   \int_0^t ds\,\psi_k\Big(\frac{s}{t}\Big)\int_{\R^3}\frac{dp}{1+p^2}
   \,|\hat{f}_j(t-s, \xi, p)|\,\psi_n\Big(\sqrt{1-(\bar{\xi}\cdot\bar{v})^2}\Big)
   \nonumber
   \\ & & \hspace{12em}\times\int_{-1}^1 d\sigma\,\frac{1}{(1+|v|\sigma)^{3/2}}\,
   \psi_m(\sqrt{1-\sigma^2}).
\end{eqnarray} 
From (\ref{qmon}) and (\ref{bern1}) we deduce that
\begin{equation}\label{berse1}
   {\|\hat{f}_j\|}_{L^2_\xi(\R^3)}
   \le C{\|f_j\|}_{L^2_x(\R^3)}\le C 2^{3j/2}{\|f_j\|}_{L^1_x(\R^3)}
   \le C 2^{3j/2}{\|f\|}_{L^1_x(\R^3)},  
\end{equation} 
and also 
\begin{equation}\label{berse2} 
   {\|\hat{f}_j\|}_{L^2_\xi(\R^3)}^2
   \le C{\|f_j\|}_{L^2_x(\R^3)}^2\le C{\|f\|}_{L^2_x(\R^3)}^2, 
\end{equation}     
where we dropped the arguments for simplicity. 

To begin with the estimate of (\ref{hatu-kmn}), 
the support of $\psi_m(\sqrt{1-\sigma^2})$ is contained in 
\[ \sigma_-=\sqrt{1-2^{-2m+2}}\le |\sigma|\le\sqrt{1-2^{-2(m+2)}}=\sigma_+. \]
Then $\sigma_+ -\sigma_-\le C2^{-2m}$ and it follows that 
\begin{eqnarray*} 
   \lefteqn{\int_{-1}^1 d\sigma\,\frac{1}{(1+|v|\sigma)^{3/2}}\,\psi_m(\sqrt{1-\sigma^2})}
   \\ & \le & \frac{2}{|v|}\Big(\frac{1}{\sqrt{1+\sigma_-|v|}}-\frac{1}{\sqrt{1+\sigma_+|v|}}
   +\frac{1}{\sqrt{1-\sigma_+|v|}}-\frac{1}{\sqrt{1-\sigma_-|v|}}\Big)
   \\ & \le & C(\sigma_+ -\sigma_-)\Big(1+(1+p^2)^{3/2}\Big)\le C2^{-2m}(1+p^2)^{3/2}. 
\end{eqnarray*} 
Thus taking $R=2^m\ge 1$,   
\begin{eqnarray*}  
   \lefteqn{\int_{\R^3}\frac{dp}{1+p^2}
   \,|\hat{f}_j(t-s, \xi, p)|\int_{-1}^1 d\sigma\,
   \frac{1}{(1+|v|\sigma)^{3/2}}\,\psi_m(\sqrt{1-\sigma^2})}
   \\ & = & \int_{|p|\le R}\frac{dp}{1+p^2}\,(\ldots)
   +\int_{|p|\ge R}\frac{dp}{1+p^2}\,(\ldots)
   \\ & \le & C2^{-2m}\int_{|p|\le R} \sqrt{1+p^2}\,|\hat{f}_j(t-s, \xi, p)|\,dp
   +C\int_{|p|\ge R}\frac{dp}{\sqrt{1+p^2}}\,\frac{1}{|v|}\,|\hat{f}_j(t-s, \xi, p)|
   \\ & \le & C2^{-2m}\int_{|p|\le R} \sqrt{1+p^2}\,|\hat{f}_j(t-s, \xi, p)|\,dp
   +CR^{-2}\int_{|p|\ge R} \sqrt{1+p^2}\,|\hat{f}_j(t-s, \xi, p)|\,dp
   \\ & \le & C2^{-2m}\int_{\R^3} \sqrt{1+p^2}\,|\hat{f}_j(t-s, \xi, p)|\,dp.
\end{eqnarray*} 
As a consequence, by (\ref{berse1}) and energy conservation (\ref{energcons}),   
\begin{eqnarray*}  
   \lefteqn{\Big\|\int_{\R^3}\frac{dp}{1+p^2}
   \,|\hat{f}_j(t-s, \cdot, p)|\int_{-1}^1 d\sigma\,\frac{1}{(1+|v|\sigma)^{3/2}}
   \,\psi_m(\sqrt{1-\sigma^2})\Big\|_{L^2_\xi(\R^3)}}
   \\ & \le & C2^{-2m}\int_{\R^3} \sqrt{1+p^2}\,{\|\hat{f}_j(t-s, \cdot, p)\|}_{L^2_\xi(\R^3)}\,dp
   \\ & \le & C2^{-2m}\,2^{3j/2}\int_{\R^3} dp\,\sqrt{1+p^2}\int_{\R^3}\,dx f(t-s, x, p)
   \\ & \le & C\,{\cal E}(0)\,2^{-2m}\,2^{3j/2}
   \\ & = & C\,2^{-2m}\,2^{3j/2}. 
\end{eqnarray*}
Using this and $\psi_n(\ldots)\le 1$ in (\ref{hatu-kmn}), we obtain 
\begin{eqnarray}\label{hatu-1}
   {\|\hat{u}_{jkmn}(t, \cdot)\|}_{L^2_\xi(\R^3)}
   & \le & C\min\Big\{1, \frac{2^{(k+m+n-j)/2}}{t^{1/2}}\Big\}
   \,2^{-2m}\,2^{3j/2}\,\int_0^t ds\,\psi_k\Big(\frac{s}{t}\Big)
   \nonumber
   \\ & \le & Ct\min\Big\{1, \frac{2^{(k+m+n-j)/2}}{t^{1/2}}\Big\}
   \,2^{-k}\,2^{-2m}\,2^{3j/2}. 
\end{eqnarray} 
Secondly, the support of $\psi_n(\sqrt{1-\tau^2})$ 
is contained in 
\[ \tau_-=\sqrt{1-2^{-2n+2}}\le |\tau|\le\sqrt{1-2^{-2(n+2)}}=\tau_+ \] 
and $\tau_+-\tau_-\le C2^{-2n}$. Thus, for $R\ge 1$, 
\begin{eqnarray*} 
   \int_{1\le |p|\le R}\frac{dp}{(1+p^2)^{3/2}}
   \,\psi_n\Big(\sqrt{1-(\bar{\xi}\cdot\bar{v})^2}\Big)
   & = & \int_{1\le |p|\le R}\frac{dp}{(1+p^2)^{3/2}}
   \,\psi_n\Big(\sqrt{1-\bar{v}_3^2}\Big)
   \\ & \le & C\int_1^R dr\,\frac{r^2}{(1+r^2)^{3/2}}\int_0^{\pi} d\varphi\,\sin\varphi
   \,\psi_n(\sqrt{1-\cos^2\varphi})
   \\ & \le & C\ln(1+R)\int_{-1}^1 \psi_n(\sqrt{1-\tau^2})\,d\tau
   \\ & \le & C\ln(1+R)\,(\tau_+-\tau_-)
   \\ & \le & C\ln(1+R)\,2^{-2n},
\end{eqnarray*} 
and similarly 
\[ \int_{|p|\le 1}\psi_n\Big(\sqrt{1-(\bar{\xi}\cdot\bar{v})^2}\Big)\,dp\le C\,2^{-2n}. \] 
Hence if we take $R=2^{n/2}\,2^{3j/4}\ge 1$, then  
\begin{eqnarray*} 
   \lefteqn{\int_{\R^3}\frac{dp}{1+p^2}
   \,|\hat{f}_j(t-s, \xi, p)|\,\psi_n\Big(\sqrt{1-(\bar{\xi}\cdot\bar{v})^2}\Big)
   \int_{-1}^1 d\sigma\,\frac{1}{(1+|v|\sigma)^{3/2}}}
   \\ & = & \int_{|p|\le 1}\frac{dp}{1+p^2}\,(\ldots)+\int_{1\le |p|\le R}\frac{dp}{1+p^2}\,(\ldots)
   +\int_{|p|\ge R}\frac{dp}{1+p^2}\,(\ldots)
   \\ & \le & C\int_{|p|\le 1} dp\,\sqrt{1+p^2}
   \,|\hat{f}_j(t-s, \xi, p)|\,\psi_n\Big(\sqrt{1-(\bar{\xi}\cdot\bar{v})^2}\Big)
   \\ & & +\,C\int_{1\le |p|\le R}\frac{dp}{\sqrt{1+p^2}}\,|\hat{f}_j(t-s, \xi, p)|
   \,\psi_n\Big(\sqrt{1-(\bar{\xi}\cdot\bar{v})^2}\Big)
   \\ & & +\,C\int_{|p|\ge R}\frac{dp}{\sqrt{1+p^2}}\,|\hat{f}_j(t-s, \xi, p)|
   \\ & \le & C\bigg(\int_{|p|\le 1}\psi_n\Big(\sqrt{1-(\bar{\xi}\cdot\bar{v})^2}\Big)\,dp\bigg)^{1/2}
   \Big(\int_{|p|\le 1} |\hat{f}_j(t-s, \xi, p)|^2\,dp\Big)^{1/2}
   \\ & & +\,C\bigg(\int_{1\le |p|\le R}\frac{dp}{(1+p^2)^{3/2}}
   \,\psi_n\Big(\sqrt{1-(\bar{\xi}\cdot\bar{v})^2}\Big)\bigg)^{1/2}
   \Big(\int_{1\le |p|\le R}\sqrt{1+p^2}\,|\hat{f}_j(t-s, \xi, p)|^2\,dp\Big)^{1/2}
   \\ & & +\,CR^{-2}\int_{|p|\ge R}\sqrt{1+p^2}\,|\hat{f}_j(t-s, \xi, p)|\,dp
   \\ & \le & C\,2^{-n}
   \Big(\int_{\R^3} |\hat{f}_j(t-s, \xi, p)|^2\,dp\Big)^{1/2}
   \\ & & +\,C(\ln(1+R))^{1/2}\,2^{-n}
   \Big(\int_{\R^3}\sqrt{1+p^2}\,|\hat{f}_j(t-s, \xi, p)|^2\,dp\Big)^{1/2}
   \\ & & +\,CR^{-2}\int_{\R^3}\sqrt{1+p^2}\,|\hat{f}_j(t-s, \xi, p)|\,dp
   \\ & \le & C(\ln(1+R))^{1/2}\,2^{-n}
   \Big(\int_{\R^3}\sqrt{1+p^2}\,|\hat{f}_j(t-s, \xi, p)|^2\,dp\Big)^{1/2}
   \\ & & +\,CR^{-2}\int_{\R^3}\sqrt{1+p^2}\,|\hat{f}_j(t-s, \xi, p)|\,dp. 
\end{eqnarray*} 
Using (\ref{berse2}), ${\|f(t)\|}_\infty\le {\|f^{(0)}\|}_\infty$ 
and (\ref{berse1}), this yields  
\begin{eqnarray*} 
   \lefteqn{\hspace{-5em}\Big\|\int_{\R^3}\frac{dp}{1+p^2}
   \,|\hat{f}_j(t-s, \cdot, p)|\,\psi_n\Big(\sqrt{1-(\bar{\xi}\cdot\bar{v})^2}\Big)
   \int_{-1}^1 d\sigma\,\frac{1}{(1+|v|\sigma)^{3/2}}\Big\|_{L^2_\xi(\R^3)}}
   \\ & \le & C(\ln(1+R))^{1/2}\,2^{-n}
   \Big(\int_{\R^3} d\xi\int_{\R^3} dp\,\sqrt{1+p^2}\,|\hat{f}_j(t-s, \xi, p)|^2\Big)^{1/2}
   \\ & & +\,CR^{-2}\int_{\R^3}\sqrt{1+p^2}\,{\|\hat{f}_j(t-s, \cdot, p)\|}_{L^2_\xi(\R^3)}\,dp
   \\ & \le & C(\ln(1+R))^{1/2}\,2^{-n}+C R^{-2}\,2^{3j/2}
   \\ & \le & C\,(\ln(1+2^{n/2}\,2^{3j/4}))^{1/2}\,2^{-n}
   \\ & \le & C\,(\sqrt{n}+\sqrt{j})\,2^{-n}. 
\end{eqnarray*} 
Due to (\ref{hatu-kmn}), and dropping $\psi_m(\ldots)\le 1$, it follows that  
\begin{eqnarray}\label{hatu-2}
   {\|\hat{u}_{jkmn}(t, \cdot)\|}_{L^2_\xi(\R^3)}
   & \le & C(0)\min\Big\{1, \frac{2^{(k+m+n-j)/2}}{t^{1/2}}\Big\}
   (\sqrt{n}+\sqrt{j})\,2^{-n}\int_0^t ds\,\psi_k\Big(\frac{s}{t}\Big) 
   \nonumber
   \\ & \le & C(0)t\min\Big\{1, \frac{2^{(k+m+n-j)/2}}{t^{1/2}}\Big\}
   (\sqrt{n}+\sqrt{j})\,2^{-k}\,2^{-n}. 
\end{eqnarray}  
Therefore if we summarize (\ref{hatu-1}) and (\ref{hatu-2}), 
we have shown (\ref{ujkmnL2}). {\hfill$\Box$}\bigskip

\begin{lemma} For $j\in\N$, 
\[ {\|u_j(t, \cdot)\|}_{L^2_x(\R^3)}
   \le C(t+\sqrt{t})\,2^{-\frac{j}{11}}. \]  
\end{lemma} 
{\bf Proof\,:} By (\ref{uhat-dec}), 
\[ {\|u_j(t, \cdot)\|}_{L^2_x(\R^3)}
   \le C{\|\hat{u}_j(t, \cdot)\|}_{L^2_\xi(\R^3)}
   \le C\sum_{k=0}^\infty\sum_{m=0}^\infty\sum_{n=0}^\infty 
   {\|\hat{u}_{jkmn}(t, \cdot)\|}_{L^2_\xi(\R^3)}. \] 
In the following, Lemma \ref{dreifrz} will be used to bound the right-hand side 
for fixed $j\in\N$ and $k\in\N_0$. Let $\alpha=\frac{16}{15}>1$ and $\eps=\frac{1}{20}\in ]0, 1[$. 
Then by Lemma \ref{dreifrz}, 
\begin{eqnarray*} 
   \sum_{m=0}^\infty\sum_{n=0}^\infty {\|\hat{u}_{jkmn}(t, \cdot)\|}_{L^2_\xi(\R^3)}
   & \le & \sum_{m,\,n} {\bf 1}_{\{m>\alpha\frac{3j}{4}\}}\,{\|\hat{u}_{jkmn}(t, \cdot)\|}_{L^2_\xi(\R^3)}
   +\sum_{m,\,n} {\bf 1}_{\{m\le\alpha\frac{3j}{4}\}}\,{\|\hat{u}_{jkmn}(t, \cdot)\|}_{L^2_\xi(\R^3)}
   \\ & \le & Ct\,2^{-k}\sum_{m=[\alpha\frac{3j}{4}]-1}^\infty\sum_{n=0}^\infty
   \,(2^{-2m}\,2^{3j/2})^{1-\eps}\,((\sqrt{n}+\sqrt{j})\,2^{-n})^\eps
   \\ & & +\,C\sqrt{t}\,2^{-k/2}\sum_{m=0}^{[\alpha\frac{3j}{4}]+1}
   \sum_{n=0}^\infty 2^{(m+n-j)/2}\,(\sqrt{n}+\sqrt{j})\,2^{-n}
   \\ & \le & Ct\,2^{-k}\,2^{3(1-\eps)j/2}\,j^{\eps/2}
   \sum_{m=[\alpha\frac{3j}{4}]-1}^\infty 2^{-2(1-\eps)m}
   \\ & & +\,C\sqrt{t}\,2^{-k/2}\,2^{-j/2}\,\sqrt{j}\sum_{m=0}^{[\alpha\frac{3j}{4}]+1} 2^{m/2}
   \\ & \le & Ct\,2^{-k}\,2^{3(1-\eps)j/2}\,j^{\eps/2}
   \,2^{-2(1-\eps)\alpha\frac{3j}{4}}
   +C\sqrt{t}\,2^{-k/2}\,2^{-j/2}\,\sqrt{j}\,2^{\,\alpha\frac{3j}{8}}
   \\ & = & Ct\,2^{-k}\,j^{\eps/2}\,2^{-(1-\eps)(\alpha-1)\frac{3j}{2}}
   +C\sqrt{t}\,2^{-k/2}\,\sqrt{j}\,2^{-\frac{j}{2}(1-\frac{3}{4}\alpha)}
   \\ & = & Ct\,2^{-k}\,j^{\frac{1}{40}}\,2^{-\frac{19}{200}\,j}
   +C\sqrt{t}\,2^{-k/2}\,\sqrt{j}\,2^{-\frac{j}{10}}
   \\ & \le & C(t+\sqrt{t})\,2^{-k/2}\,2^{-\frac{j}{11}}. 
\end{eqnarray*} 
Summation on $k\in\N_0$ concludes the proof of the lemma. 
{\hfill$\Box$}\bigskip

Now we are in position to finish the proof of Theorem \ref{mainthm}. 
To summarize, we have seen that 
\[ |E_T(t, x)|\le C u(t, x),\quad u=\sum_{j=0}^\infty u_j,
   \quad {\|u_j(t, \cdot)\|}_{L^2_x}\le C(t+\sqrt{t})\,2^{-\frac{j}{11}} \] 
for $j\in\N$. Clearly one also has ${\|u_0(t, \cdot)\|}_{L^2_x}\le Ct$. 
Let $H^s_x(\R^3)$ denote the standard (inhomogeneous) $L^2_x$-based Sobolev space 
of order $s$. Then by the inhomogeneous Sobolev embedding theorem 
and by Plancherel's theorem, for $2<q<\infty$, $s>0$ 
and $\frac{1}{2}\le\frac{1}{q}+\frac{s}{3}$: 
\begin{eqnarray*} 
   {\|E_T(t, \cdot)\|}_{L^q_x} & \le & C {\|u(t, \cdot)\|}_{L^q_x}
   \\ & \le & C{\|u(t, \cdot)\|}_{H^s_x}
   \\ & \le & C\Big[{\|u_0(t, \cdot)\|}_{L^2_x}
   +\Big(\sum_{j=1}^\infty 2^{2sj}\|u_j(t, \cdot)\|_{L^2_x}^2\Big)^{1/2}\Big]
   \\ & \le & C(t)\Big[1
   +\Big(\sum_{j=1}^\infty 2^{2j(s-\frac{1}{11})}\Big)^{1/2}\Big]
   \\ & \le & C(t),
\end{eqnarray*} 
provided that $s<\frac{1}{11}$. Hence $q=2+\delta$ is possible, 
and for instance $s=\frac{1}{12}$ and $\delta=\frac{2}{17}$ 
is a suitable choice. {\hfill$\Box$}\bigskip 


\setcounter{equation}{0}

\section{Appendix: Explicit form of the kernels}
\label{kersect} 

To make this paper self-contained, we will include the following formulas; 
see \cite[Section II]{glstr} and \cite[(A13), (A14), (A3)]{schaeffer:86}. 
The fields $E$ and $B$ can be written as 
\begin{eqnarray*}
   E & = & E_D+E_{DT}+E_T+E_S,
   \\ B & = & B_D+B_{DT}+B_T+B_S,
\end{eqnarray*}
where
\begin{eqnarray*}
   E_D(t, x) & = &
   \partial_t\bigg(\frac{t}{4\pi}\int_{|\omega|=1} E^{(0)}(x+t\omega)\,d\omega\bigg)
   \\ & & +\,\frac{t}{4\pi}\int_{|\omega|=1}\partial_t E(0, x+t\omega)\,d\omega, \\[1ex]
   E_{DT}(t, x) & = & -\frac{1}{t}\int_{|y|=t}\int_{\R^3} K_{E,\,DT}(\omega, v)
   f^{(0)}(x+y, p)\,dp\,d\sigma(y), \\[1ex] 
   B_D(t, x) & = &
   \partial_t\bigg(\frac{t}{4\pi}\int_{|\omega|=1} B^{(0)}(x+t\omega)\,d\omega\bigg)
   \\ & & +\,\frac{t}{4\pi}\int_{|\omega|=1}\partial_t B(0, x+t\omega)\,d\omega, \\[1ex]
   B_{DT}(t, x) & = & \frac{1}{t}\int_{|y|=t}\int_{\R^3} K_{B,\,DT}(\omega, v)
   f^{(0)}(x+y, p)\,dp\,d\sigma(y),  
\end{eqnarray*} 
are the data terms. In addition, 
\begin{eqnarray*} 
   E_T(t, x) & = & -\int_{|y|\le t}\frac{dy}{|y|^2}
   \int_{\R^3} dp\,K_{E,\,T}(\omega, v) f(t-|y|, x+y, p), \\[1ex] 
   E_S(t, x) & = & -\int_{|y|\le t}\frac{dy}{|y|}
   \int_{\R^3} dp\,K_{E,\,S}(\omega, v)\,(Lf)(t-|y|, x+y, p),
\end{eqnarray*}
and
\begin{eqnarray*}
   B_T(t, x) & = & \int_{|y|\le t}\frac{dy}{|y|^2}
   \int_{\R^3} dp\,K_{B,\,T}(\omega, v) f(t-|y|, x+y, p), \\[1ex] 
   B_S(t, x) & = & \int_{|y|\le t}\frac{dy}{|y|}
   \int_{\R^3} dp\,K_{B,\,S}(\omega, v)\,(Lf)(t-|y|, x+y, p),
\end{eqnarray*}
defining $\omega=|y|^{-1}y$ and $L=E+v\wedge B$. 
The kernels are 
\begin{eqnarray*}
   K_{E,\,DT}(\omega, v) & = & (1+v\cdot\omega)^{-1}
   (\omega-(v\cdot\omega)v), \\
   K_{E,\,T}(\omega, v)  & = & (1+p^2)^{-1}(1+v\cdot\omega)^{-2}(v+\omega), \\ 
   K_{E,\,S}(\omega, v)  & = & (1+p^2)^{-1/2}(1+v\cdot\omega)^{-2}
   \\ & & \Big[(1+v\cdot\omega)+((v\cdot\omega)\omega-v)\otimes v
   -(v+\omega)\otimes\omega\Big]\in\R^{3\times 3},
\end{eqnarray*}
and 
\begin{eqnarray*}
   K_{B,\,DT}(\omega, v) & = & -(1+v\cdot\omega)^{-1}(v\wedge\omega), \\ 
   K_{B,\,T}(\omega, v) & = & -(1+p^2)^{-1}(1+v\cdot\omega)^{-2}(v\wedge\omega), \\
   K_{B,\,S}(\omega, v) & = & (1+p^2)^{-1/2}(1+v\cdot\omega)^{-2}
   \\ & & \Big[(1+v\cdot\omega)\,\omega\wedge (\ldots)
   -(v\wedge\omega)\otimes (v+\omega)\Big]\in\R^{3\times 3}.
\end{eqnarray*}
\medskip 

\noindent 
{\bf Proof of (\ref{beh1}) and (\ref{beh2})\,:} The bound (\ref{beh1}) 
is immediate from 
\[ |v+\omega|=(v^2+2(v\cdot\omega)+1)^{1/2}
   \le\sqrt{2}\,(1+v\cdot\omega)^{1/2}. \] 
Regarding (\ref{beh2}), we use that 
\begin{eqnarray*}
   \Big[((v\cdot\omega)\omega-v)\otimes v-(v+\omega)\otimes\omega\Big]z
   & = & (v\cdot z)((v\cdot\omega)\omega-v)-(\omega\cdot z)(v+\omega)
   \\ & = & -\,(\omega-(v\cdot\omega)v)\cdot z\,(v+\omega)
   -(1+v\cdot\omega)(v\cdot z)\,v
\end{eqnarray*}
and 
\begin{eqnarray*}
   |\omega-(v\cdot\omega)v| & = & (1-2(v\cdot\omega)^2+(v\cdot\omega)^2 v^2)^{1/2}
   \\ & \le & (1-(v\cdot\omega)^2)^{1/2}\le\sqrt{2}\,(1+v\cdot\omega)^{1/2}. 
\end{eqnarray*} 
This yields the claim.
{\hfill$\Box$}\bigskip 


\end{document}